\documentclass[letterpaper]{amsproc}
\usepackage{amssymb}
\usepackage{amsfonts}

\theoremstyle{plain}

\newtheorem{corollary}{Corollary}

\newtheorem{definition}{Definition}

\newtheorem{lemma}{Lemma}

\newtheorem{proposition}{Proposition}

\newtheorem{theorem}{Theorem}
\numberwithin{equation}{section}
\input{tcilatex}

\begin{document}
\title[Ergodic Theorems for Lower Probabilities]{Ergodic Theorems for Lower
Probabilities}
\author{S.Cerreia--Vioglio, F. Maccheroni, and M. Marinacci}
\address{Universit\`{a} Bocconi}
\thanks{Corresponding Author: Fabio Maccheroni \TEXTsymbol{<}%
fabio.maccheroni@unibocconi.it\TEXTsymbol{>}, U. Bocconi, via Sarfatti 25,
20136, Milano, ITALY. The authors\ gratefully acknowledge the financial
support of MIUR (PRIN grant 20103S5RN3\_005).}

\begin{abstract}
We establish an Ergodic Theorem for lower probabilities, a generalization of
standard probabilities widely used in applications. As a by-product, we
provide a version for lower probabilities of the Strong Law of Large Numbers.
\end{abstract}

\maketitle

\section{Introduction}

The purpose of this paper is to state and prove an Ergodic Theorem for lower
probabilities: a class of monotone set functions that are not necessarily
additive and are\ widely used in applications where standard additive
probabilities turn out to be inadequate (for applications in Economics see
Marinacci and Montrucchio \cite{MM}, for applications in Statistics see
Walley \cite{Wal}).

We consider a measurable space $\left( \Omega,\mathcal{F}\right) $, endowed
with an $\mathcal{F}\backslash \mathcal{F}$-measurable transformation $%
\tau:\Omega \rightarrow \Omega$, and a (continuous) lower probability $\nu:%
\mathcal{F}\rightarrow \left[ 0,1\right] $. We study four different notions
of invariance for lower probabilities (Definitions \ref{def:cap-inv}-\ref%
{def:cap-bay}). They are equivalent in the additive case, and so are genuine
generalizations to the nonadditive setting of the usual concept of
invariance.

The most natural definition of invariance for a lower probability $\nu $
(Definition \ref{def:cap-inv}) requires that%
\begin{equation*}
\nu \left( A\right) =\nu \left( \tau ^{-1}\left( A\right) \right) \qquad
\forall A\in \mathcal{F}.
\end{equation*}%
It is the weakest form of invariance for the nonadditive case. Nevertheless,
it is still possible to derive a version of the Ergodic Theorem (Theorem \ref%
{thm:w-inv-erg}). In other words, if $\nu $ is an invariant lower
probability, then for each real valued, bounded, and measurable function $%
f:\Omega \rightarrow 
\mathbb{R}
$ the limit 
\begin{equation*}
\lim_{n}\frac{1}{n}\sum_{k=1}^{n}f\circ \tau ^{k-1}\left( \omega \right)
\end{equation*}%
exists on a set that has measure $1$ with respect to $\nu $. If, in
addition, $\nu $ is ergodic, we are able to provide bounds for such limit in
terms of lower\ and upper Choquet integrals.

Under the stronger notions of invariance (Definitions \ref{def:cap-str-inv}-%
\ref{def:cap-bay}), the previous result can be strengthened in several ways.
First, we develop a nonadditive version of Kingman's super-subadditive
ergodic theorem (Theorem \ref{thm:sup-erg}). Second, when $\left( \Omega ,%
\mathcal{F}\right) $ is a standard measurable space we can better
characterize the limit of time averages (Corollary \ref{cor:erg}).

As an application of our main result, we establish a nonadditive version of
the Strong Law of Large Numbers (Theorem \ref{thm:SLLN}) for stationary and
ergodic processes.

\section{Mathematical Preliminaries}

\subsection{Set functions}

Consider a measurable space $\left( S,\Sigma \right) $, where $S$ is a
nonempty set and $\Sigma $ is a $\sigma $-algebra of subsets of $S$. Subsets
of $S$ are understood to be in $\Sigma $ even when\ not stated explicitly. A
set function $\nu :\Sigma \rightarrow \left[ 0,1\right] $ is

\begin{enumerate}
\item[(i)] \emph{a capacity} if $\nu \left( \emptyset \right) =0$, $\nu
\left( S\right) =1$, and $\nu \left( A\right) \leq \nu \left( B\right) $ for
all $A$ and $B\ $such that $A\subseteq B$;

\item[(ii)] \emph{convex }if $\nu \left( A\cup B\right) +\nu \left( A\cap
B\right) \geq \nu \left( A\right) +\nu \left( B\right) $ for all $A$ and $B$;

\item[(iii)] \emph{additive }if $\nu \left( A\cup B\right) =\nu \left(
A\right) +\nu \left( B\right) $ for all disjoint $A$ and $B$;

\item[(iv)] \emph{continuous }if $\lim_{n\rightarrow \infty}\nu \left(
A_{n}\right) =\nu \left( A\right) $ whenever either $A_{n}\downarrow A$ or $%
A_{n}\uparrow A$;

\item[(v)] \emph{continuous} \emph{at} $S$ if $\lim_{n\rightarrow \infty
}\nu \left( A_{n}\right) =\nu \left( S\right) $ whenever $A_{n}\uparrow S$;

\item[(vi)] \emph{a probability\ }if it is an additive capacity;

\item[(vii)] \emph{a probability measure} if it is a probability which is
continuous at $S$.
\end{enumerate}

We denote by $\Delta \left( S,\Sigma \right) $ the set of all probabilities
on $\Sigma $ and by $\Delta ^{\sigma }\left( S,\Sigma \right) $ the set of
all probability measures on $\Sigma $. We endow both sets with the relative
topology induced by the weak* topology.\footnote{%
Recall that a net $\left\{ P_{\alpha }\right\} _{\alpha \in I}$ converges to 
$P$, in the weak* topology, if and only if $P_{\alpha }\left( A\right)
\rightarrow P\left( A\right) $ for all $A\in \Sigma $. The weak* topology is
thus the restriction to $\Delta \left( S,\Sigma \right) $ of the topology $%
\sigma \left( ba\left( S,\Sigma \right) ,B\left( S,\Sigma \right) \right) $
where $B\left( S,\Sigma \right) $ is the space of all real valued, bounded,
and $\Sigma $\textit{-measurable} functions on $S$ and $ba\left( S,\Sigma
\right) $ is the set of all bounded and finitely additive set functions on $%
\Sigma $. In the case of $S$ being a Polish space and $\Sigma $ the Borel $%
\sigma $-algebra, the above topology should not be confused with the
topology generated by real valued, bounded, and \textit{continuous }%
functions on $S$.} Given $\mathcal{M}\subseteq \Delta ^{\sigma }\left(
S,\Sigma \right) $, we assume that $\mathcal{M}$ is endowed with the $\sigma 
$-algebra $\mathcal{A}_{\mathcal{M}}$ which is the smallest $\sigma $%
-algebra that makes the evaluations $P\mapsto P\left( A\right) $ measurable
for all $A\in \Sigma $. A set function $\nu :\Sigma \rightarrow \left[ 0,1%
\right] $ is

\begin{enumerate}
\item[(viii)] \emph{a lower probability (measure)} if there exists a compact
set $\mathcal{M}\subseteq \Delta^{\sigma}\left( S,\Sigma \right) $ such that%
\begin{equation*}
\nu \left( A\right) =\min_{P\in \mathcal{M}}P\left( A\right) \qquad \forall
A\in \Sigma.
\end{equation*}
\smallskip
\end{enumerate}

Given a capacity $\nu $, its conjugate $\bar{\nu}:\Sigma \rightarrow \left[
0,1\right] $ is given by%
\begin{equation*}
\bar{\nu}\left( A\right) =1-\nu \left( A^{c}\right) \qquad \forall A\in
\Sigma .
\end{equation*}%
It is immediate to verify that if $\nu $ is a lower probability, then%
\begin{equation}
\bar{\nu}\left( A\right) =\max_{P\in \mathcal{M}}P\left( A\right) \qquad
\forall A\in \Sigma .  \label{eq:conj}
\end{equation}%
The core of a capacity $\nu $ is the weak* compact set defined by%
\begin{equation*}
\mathrm{core}\left( \nu \right) =\left\{ P\in \Delta \left( S,\Sigma \right)
:P\geq \nu \right\} ,
\end{equation*}%
that is, the core is the collection of all probabilities that setwise
dominate $\nu $. A capacity\ $\nu :\Sigma \rightarrow \left[ 0,1\right] $ is

\begin{enumerate}
\item[(ix)] \emph{exact} if $\mathrm{core}\left( \nu \right) \not =\emptyset$
and $\nu \left( A\right) =\min_{P\in \mathrm{core}\left( \nu \right)
}P\left( A\right) $ for each $A$.
\end{enumerate}

If $\nu $ is a convex capacity continuous at $S$, then $\nu $ is exact and $%
\emptyset \not=\mathrm{core}\left( \nu \right) \subseteq \Delta ^{\sigma
}\left( S,\Sigma \right) $ (see \cite[Lemma 2 and Theorem 1]{Del}, \cite[%
Theorem 3.2]{Sch}, and \cite[Theorem 4.2 and Theorem 4.7]{MM}). In
particular, $\nu $ is a lower probability where $\mathcal{M}=\mathrm{core}%
\left( \nu \right) $. Conversely, if $\nu $ is a lower probability, then $%
\nu $ is exact, continuous at $S$, and $\mathcal{M}\subseteq \mathrm{core}%
\left( \nu \right) \subseteq \Delta ^{\sigma }\left( S,\Sigma \right) $.
Nevertheless, being exact does not automatically imply being convex. An
exact capacity continuous at $S$ is continuous. Finally, we say that a
statement about a random element holds $\nu -a.s.$ if and only if there
exists an event $A$ such that $\nu \left( A\right) =1$ and the statement
holds for all $s\in A$.

\subsection{Integrals}

We denote by $B\left( S,\Sigma \right) $ the set of all bounded and $\Sigma $%
-measurable functions from $S$ to $%
\mathbb{R}
$. A capacity $\nu $ induces a functional on $B\left( S,\Sigma \right) $ via
the Choquet integral, defined for all $f\in B\left( S,\Sigma \right) $ by:%
\begin{equation*}
\int_{S}fd\nu =\int_{0}^{\infty }\nu \left( \left\{ s\in S:f\left( s\right)
\geq t\right\} \right) dt+\int_{-\infty }^{0}\left[ \nu \left( \left\{ s\in
S:f\left( s\right) \geq t\right\} \right) -\nu \left( S\right) \right] dt
\end{equation*}%
where the right hand side integrals are (improper) Riemann integrals. If $%
\nu $ is additive, then the Choquet integral reduces to the standard
additive integral. It is also routine to check that $-\int_{S}fd\nu
=\int_{S}-fd\bar{\nu}$ for all $f\in B\left( S,\Sigma \right) $. It is well
known (see \cite[Lemma 2]{Del}, \cite[Proposition 3]{Sch1}, and \cite[%
Theorem 4.7]{MM}) that if $\nu $ is a convex capacity, then%
\begin{equation*}
\int_{S}fd\nu =\min_{P\in \mathrm{core}\left( \nu \right) }\int_{S}fdP\text{
and }\int_{S}fd\bar{\nu}=\max_{P\in \mathrm{core}\left( \nu \right)
}\int_{S}fdP\qquad \forall f\in B\left( S,\Sigma \right) .
\end{equation*}

In the rest of the paper, we consider three measurable spaces $\left(
S,\Sigma \right) $. The first one is $\left( \Omega ,\mathcal{F}\right) $
which we interpret as the space where ultimately uncertainty lives. Given a
set $\mathcal{P}\subseteq \Delta ^{\sigma }\left( \Omega ,\mathcal{F}\right) 
$, the second space will be $\left( \mathcal{P},\mathcal{A}_{\mathcal{P}%
}\right) $ which we interpret as the space of all possible probability
models equipped with the $\sigma $-algebra $\mathcal{A}_{\mathcal{P}}$\
discussed above. Finally, given a real valued and $\mathcal{F}$-measurable
stochastic process $\left \{ f_{n}\right \} _{n\in 
\mathbb{N}
}$ on $\Omega $, we will consider the space $\left( 
\mathbb{R}
^{%
\mathbb{N}
},\sigma \left( \mathcal{C}\right) \right) $, which we will interpret as the
space of observations endowed with the $\sigma $-algebra generated by the
algebra of cylinders $\mathcal{C}$.

\subsection{Prior and Predictive Capacities}

Given a set $\mathcal{P}\subseteq \Delta^{\sigma}\left( \Omega,\mathcal{F}%
\right) $, a prior is a capacity $\rho:\mathcal{A}_{\mathcal{P}}\rightarrow%
\left[ 0,1\right] $. The associated predictive is the capacity $\nu_{\rho}:%
\mathcal{F}\rightarrow \left[ 0,1\right] $ defined by%
\begin{equation*}
\nu_{\rho}\left( A\right) =\int_{\mathcal{P}}P\left( A\right) d\rho \left(
P\right) \qquad \forall A\in \mathcal{F}.
\end{equation*}
If $\rho$ is additive and continuous at $\mathcal{P}$, then $\rho$ is a
prior and $\nu_{\rho}$ is a predictive in the traditional sense. We denote
capacities that are additive and continuous at $\mathcal{P}$ by $\pi$. Given
a set $\mathcal{P}$, we denote the set of strong extreme points of $\mathcal{%
P}$ by $\mathcal{S}\left( \mathcal{P}\right) $.\footnote{%
Recall that $P\in \mathcal{P}$ is a strong extreme point of $\mathcal{P}$ if
and only if the Dirac at $P$ (i.e., $\delta_{P}$) is the only probability
measure $\pi:\mathcal{A}_{\mathcal{P}}\rightarrow \left[ 0,1\right] $ such
that $P\left( A\right) =\int_{\mathcal{P}}Q\left( A\right) d\pi \left(
Q\right) $ for each $A\in \mathcal{F}$.}

\section{Ergodic Theorems\label{sec:erg}}

\subsection{Invariant Capacities}

In this section, we consider a measurable space $\left( \Omega ,\mathcal{F}%
\right) $. We also consider a transformation $\tau :\Omega \rightarrow
\Omega $ which is $\mathcal{F}/\mathcal{F}$-measurable. Recall that a
probability measure $P$ is ($\tau $-)invariant if and only if%
\begin{equation}
P\left( A\right) =P\left( \tau ^{-1}\left( A\right) \right) \qquad \forall
A\in \mathcal{F}.  \label{eq:inv}
\end{equation}%
We denote by $\mathcal{I}$ the set of all probability measures that satisfy (%
\ref{eq:inv}) and by $\mathcal{G}$ the set of all invariant events of $%
\mathcal{F}$, that is, $A\in \mathcal{G}$ if and only if $A\in \mathcal{F\ }$%
and $\tau ^{-1}\left( A\right) =A$. An invariant probability measure $P$ is
said to be ergodic if and only if $P\left( \mathcal{G}\right) =\left\{
0,1\right\} $. Similarly, we say that a capacity $\nu $ is ergodic if and
only if $\nu \left( \mathcal{G}\right) =\left\{ 0,1\right\} $. We denote by $%
\mathcal{S}\left( \mathcal{I}\right) $ the subset of $\mathcal{I}$ such that%
\begin{equation*}
\mathcal{S}\left( \mathcal{I}\right) =\left\{ P\in \mathcal{I}:P\left( 
\mathcal{G}\right) =\left\{ 0,1\right\} \right\} .
\end{equation*}%
If $\left( \Omega ,\mathcal{F}\right) $ is a standard measurable space, then
it can be checked that $\mathcal{S}\left( \mathcal{I}\right) $ is the set of
strong extreme points of $\mathcal{I}$ (see Dynkin \cite{Dyn}). Finally,
following Dunford and Schwartz \cite[pp. 723-724]{DS} (see also Dowker \cite%
{Dow}), we say that a probability measure $P$ is potentially ($\tau $%
-)invariant if and only if there exists a probability measure $\hat{P}\in 
\mathcal{I}$ such that%
\begin{equation*}
P\left( E\right) =\hat{P}\left( E\right) \qquad \forall E\in \mathcal{G}.
\end{equation*}%
We denote the set of potentially invariant probability measures by $\mathcal{%
PI}$.

\bigskip

Next, we propose four notions of ($\tau$-)invariance for a capacity.

\begin{definition}
\label{def:cap-inv}A capacity $\nu$ is invariant if and only if for each $%
A\in \mathcal{F}$%
\begin{equation*}
\nu \left( A\right) =\nu \left( \tau^{-1}\left( A\right) \right) .
\end{equation*}
\end{definition}

\begin{definition}
\label{def:cap-str-inv}A capacity $\nu$ is strongly invariant if and only if
for each $A\in \mathcal{F}$%
\begin{equation*}
\nu \left( A\backslash \tau^{-1}\left( A\right) \right) =\bar{\nu}\left(
\tau^{-1}\left( A\right) \backslash A\right) \text{\ and\ }\nu \left(
\tau^{-1}\left( A\right) \backslash A\right) =\bar{\nu}\left( A\backslash
\tau^{-1}\left( A\right) \right) .
\end{equation*}
\end{definition}

\begin{definition}
\label{def:cap-fun-inv}A lower probability $\nu$ is functionally invariant
if and only if $\mathcal{M}\subseteq \mathcal{I}$.
\end{definition}

The fourth definition also describes a procedure in which invariant
capacities can be constructed. Such a procedure is a robust Bayesian
procedure (see Berger \cite{Ber} and Shafer \cite{Sha}).

\begin{definition}
\label{def:cap-bay}A capacity $\nu$ is robustly invariant if and only if $%
\nu=\nu_{\rho}$ for some convex capacity $\rho:\mathcal{A}_{\mathcal{S}%
\left( \mathcal{I}\right) }\rightarrow \left[ 0,1\right] $.
\end{definition}

It can be shown that if $\left( \Omega,\mathcal{F}\right) $ is a standard
measurable space and $\nu$ is robustly invariant and continuous at\ $\Omega$%
, then it is a lower probability. In the next two results, we will clarify
the connection between these four notions of invariance.

\begin{proposition}
\label{pro:inv-low}Let $\left( \Omega,\mathcal{F}\right) $ be a standard
measurable space and $\nu$ a lower probability. The following statements are
true:

\begin{enumerate}
\item If $\nu$ is strongly invariant, then $\nu$ is functionally
invariant\thinspace and $\mathrm{core}\left( \nu \right) \subseteq \mathcal{I%
} $.

\item If $\nu$ is robustly invariant, then $\nu$ is functionally invariant.

\item If $\nu$ is functionally invariant and $\mathcal{M}\in \mathcal{A}_{%
\mathcal{S}\left( \mathcal{I}\right) }$, then $\nu$ is robustly invariant
and ergodic.

\item If $\nu$ is functionally invariant, then $\nu$ is invariant.
\end{enumerate}
\end{proposition}

The connection among some of these notions of invariance becomes sharper
when $\nu$ is convex.

\begin{theorem}
\label{thm:inv-low-con}Let $\left( \Omega,\mathcal{F}\right) $ be a standard
measurable space and $\nu$ a convex capacity continuous at $\Omega$. The
following statements are equivalent:

\begin{enumerate}
\item[(i)] $\nu$ is strongly invariant;

\item[(ii)] $\nu$ is functionally invariant\thinspace and $\mathrm{core}%
\left( \nu \right) \subseteq \mathcal{I}$;

\item[(iii)] $\nu$ robustly invariant and $\mathrm{core}\left( \nu \right)
\subseteq \mathcal{I}$;

\item[(iv)] $\mathrm{core}\left( \nu \right) \subseteq \mathcal{I}$.
\end{enumerate}
\end{theorem}

As a corollary, we obtain that the four definitions coincide with the usual
definition of invariance when $\nu $ is a probability measure. Under
additional assumptions on $\Omega $ and $\tau $, in the additive case, the
equivalence between points (i) and (iii) follows by an application of the
Choquet-Bishop-de Leeuw theorem (see Phelps \cite{Phe}). In our case, the
equivalence between points (i) and (iii) could be proven by developing a
nonadditive version of the Choquet-Bishop-de Leeuw theorem. This can be
achieved\ by using the techniques contained in Cerreia-Vioglio, Maccheroni,
Marinacci, and Montrucchio \cite{CMMM5}.\ Finally, in the next section, we
show that, if $\nu $ is an invariant lower probability, then its core$\ $%
must be contained in $\mathcal{PI}$.

\subsection{Ergodic Theorem}

Given the notions of invariance previously discussed, we could then ask
ourselves if suitable ergodic theorems can be developed for nonadditive
probabilities. In light of Proposition \ref{pro:inv-low} and Theorem \ref%
{thm:inv-low-con}, an immediate dichotomy presents. In fact, the notion of
invariance of Definition \ref{def:cap-inv} stands separate from, and it is
actually weaker than, the other notions of strong, robust, and functional
invariance, even in the convex case. Theorem \ref{thm:w-inv-erg}\ only
assumes the weak form of invariance of Definition \ref{def:cap-inv}. On the
other hand, Corollary \ref{cor:erg} assumes strong invariance. Strong
invariance, paired with the convexity of $\nu $ and $\left( \Omega ,\mathcal{%
F}\right) $ being standard, allows us to provide a sharper version of
Theorem \ref{thm:w-inv-erg}.

\begin{theorem}
\label{thm:w-inv-erg}Let $\left( \Omega,\mathcal{F}\right) $ be a measurable
space and $\nu$ a lower probability. If $\nu$ is invariant, then for each $%
f\in B\left( \Omega,\mathcal{F}\right) $ there exists $f^{\star}\in B\left(
\Omega,\mathcal{G}\right) $ such that%
\begin{equation*}
\lim_{n}\frac{1}{n}\dsum \limits_{k=1}^{n}f\left( \tau^{k-1}\left( \omega
\right) \right) =f^{\star}\left( \omega \right) \qquad \nu-a.s.
\end{equation*}
Moreover, if $\nu$ is ergodic, then%
\begin{equation*}
\nu \left( \left \{ \omega \in \Omega:\int_{\Omega}f^{\star}d\nu \leq \lim
_{n}\frac{1}{n}\dsum \limits_{k=1}^{n}f\left( \tau^{k-1}\left( \omega
\right) \right) \leq \int_{\Omega}f^{\star }d\bar{\nu}\right \} \right) =1.
\end{equation*}
\end{theorem}

As a corollary, we are able to show a necessary property that $\mathrm{core}%
\left( \nu \right) $ of an invariant lower probability $\nu $ must satisfy
(cf. Proposition \ref{pro:inv-low}). Clearly, it is not a characterization
since it is well known that there are probability measures that are
potentially invariant, but not invariant.

\begin{corollary}
\label{cor:w-inv-erg}If a lower probability $\nu $ is invariant, then $%
\mathrm{core}\left( \nu \right) \subseteq \mathcal{PI}$.
\end{corollary}

As a second corollary, we discuss the ergodic theorem for convex and
strongly invariant capacities. Compared to Theorem \ref{thm:w-inv-erg}, the
following corollary assumes $\nu $ convex and a stronger form of invariance
that, in turn, yield\ a limit function $f^{\star }$ which has more
properties. These properties naturally generalize the ones found in the
Individual Ergodic Theorem of Birkhoff. In this case, convergence of
empirical averages is a simple consequence of Birkhoff's theorem applied to
each probability in $\mathrm{core}\left( \nu \right) $. Nevertheless, the
relation between $f$ and $f^{\star }$ in terms of Choquet expectations is
not immediate at first sight.\ A similar comment applies to Theorem \ref%
{thm:sup-erg}.

\begin{corollary}
\label{cor:erg}Let $\left( \Omega,\mathcal{F}\right) \,$be a standard
measurable space and $\nu$ a convex capacity continuous at $\Omega$. If $\nu$
is strongly invariant, then for each $f\in B\left( \Omega,\mathcal{F}\right) 
$ there exists $f^{\star}\in B\left( \Omega,\mathcal{G}\right) $ such that%
\begin{equation}
\lim_{n}\frac{1}{n}\dsum \limits_{k=1}^{n}f\left( \tau^{k-1}\left( \omega
\right) \right) =f^{\star}\left( \omega \right) \text{\qquad}\nu-a.s.
\label{eq:erg}
\end{equation}
\noindent Moreover,

\begin{enumerate}
\item For each $P\in \mathcal{I}$, $f^{\star}$ is a version of the
conditional expectation of $f$ given $\mathcal{G}$.

\item $\int_{\Omega}f^{\star}d\nu=\int_{\Omega}fd\nu$.

\item If $\nu $ is ergodic, then%
\begin{equation*}
\nu \left( \left \{ \omega \in \Omega :\int_{\Omega }fd\nu \leq \lim_{n}%
\frac{1}{n}\dsum \limits_{k=1}^{n}f\left( \tau ^{k-1}\left( \omega \right)
\right) \leq \int_{\Omega }fd\bar{\nu}\right \} \right) =1.
\end{equation*}
\end{enumerate}
\end{corollary}

\subsection{Subadditive Ergodic Theorem\label{sub:sub-sup}}

Next we turn to a Subadditive Ergodic Theorem for lower probabilities.

\begin{definition}
A sequence $\left \{ S_{n}\right \} _{n\in%
\mathbb{N}
}$ of $\mathcal{F}$-measurable random variables is superadditive (resp.,
subadditive) if and only if%
\begin{equation*}
S_{n+k}\geq S_{n}+S_{k}\circ \tau^{n}\text{ (resp., }\leq \text{)}\qquad
\forall n,k\in%
\mathbb{N}
.
\end{equation*}
The sequence $\left \{ S_{n}\right \} _{n\in%
\mathbb{N}
}$ is additive if and only if it is superadditive and subadditive.
\end{definition}

Consider an $\mathcal{F}$-measurable function $f:\Omega \rightarrow 
\mathbb{R}
$. If we define $\left \{ S_{n}\right \} _{n\in 
\mathbb{N}
}$ by%
\begin{equation}
S_{n}=\dsum \limits_{k=1}^{n}f\circ \tau ^{k-1}\qquad \forall n\in 
\mathbb{N}
,  \label{eq:add}
\end{equation}%
then we have that $\left \{ S_{n}\right \} _{n\in 
\mathbb{N}
}$ is an additive sequence. The opposite is also true, that is, if $\left \{
S_{n}\right \} _{n\in 
\mathbb{N}
}$ is additive, then it takes the form (\ref{eq:add}) for some $\mathcal{F}$%
-measurable real valued function $f$. On the other hand, if we take $%
\left
\{ S_{n}\right \} _{n\in 
\mathbb{N}
}$ as in (\ref{eq:add}) and we consider $\left \{ \left \vert
S_{n}\right
\vert \right \} _{n\in 
\mathbb{N}
}\,$we obtain a genuine subadditive sequence. Note that if $f\in B\left(
\Omega ,\mathcal{F}\right) $, then we also have that there exists $\lambda
\in 
\mathbb{R}
$ such that%
\begin{equation}
-\lambda n\leq S_{n}\left( \omega \right) \leq \lambda n\qquad \forall
\omega \in \Omega .  \label{eq:bdd}
\end{equation}%
Similarly, we have that $-\lambda n\leq \left \vert S_{n}\right \vert \leq
\lambda n$ for all $n\in 
\mathbb{N}
$.

\begin{theorem}
\label{thm:sup-erg}Let $\left( \Omega,\mathcal{F}\right) \,$be a standard
measurable space and $\nu$ a lower probability. If $\left \{ S_{n}\right \}
_{n\in%
\mathbb{N}
}$ is either a superadditive or a subadditive sequence that satisfies (\ref%
{eq:bdd}) and if $\nu$ is functionally invariant, then there exists $%
f^{\star}\in B\left( \Omega,\mathcal{G}\right) $ such that%
\begin{equation*}
\lim_{n}\frac{S_{n}}{n}=f^{\star}\text{\qquad}\nu-a.s.
\end{equation*}
\noindent Moreover,

\begin{enumerate}
\item If $\nu$ is convex and strongly invariant and $\left \{ S_{n}\right \}
_{n\in%
\mathbb{N}
}$ superadditive, then $\int_{\Omega}f^{\star}d\nu=\sup_{n\in%
\mathbb{N}
}\int_{\Omega}\frac{S_{n}}{n}d\nu$.

\item If $\nu$ is convex and strongly invariant and $\left \{ S_{n}\right \}
_{n\in%
\mathbb{N}
}$ subadditive, then $\int_{\Omega}f^{\star}d\bar{\nu}=\inf_{n}\int_{\Omega }%
\frac{S_{n}}{n}d\bar{\nu}$.

\item If $\nu $ is ergodic and $\left \{ S_{n}\right \} _{n\in 
\mathbb{N}
}$ is either subadditive or superadditive, then%
\begin{equation*}
\nu \left( \left \{ \omega \in \Omega :\int_{\Omega }f^{\star }d\nu \leq
\lim_{n}\frac{S_{n}\left( \omega \right) }{n}\leq \int_{\Omega }f^{\star }d%
\bar{\nu}\right \} \right) =1.
\end{equation*}
\end{enumerate}
\end{theorem}

\section{Strong Law of Large Numbers\label{sec:SLLN}}

As an application of Theorem \ref{thm:w-inv-erg}, we provide a nonadditive
version of the Strong Law of Large Numbers. Before doing so, we need to
introduce some notation and terminology. Consider a sequence of real valued,
bounded, and measurable random variables $\mathbf{f}=\left \{ f_{n}\right \}
_{n\in%
\mathbb{N}
}\subseteq B\left( \Omega,\mathcal{F}\right) $. We denote by $\mathcal{T}$
the tail $\sigma$-algebra $\dbigcap \limits_{k\in \mathbb{N} }\sigma \left(
f_{k},f_{k+1},...\right) $.

\begin{definition}
Given a capacity $\nu$, we say that $\mathbf{f}=\left \{ f_{n}\right \}
_{n\in%
\mathbb{N}
}$ is stationary if and only if for each $n\in%
\mathbb{N}
$, for each $k\in%
\mathbb{N}
_{0}$, and for each Borel subset $B$ of $%
\mathbb{R}
^{k+1}$%
\begin{equation}
\nu \left( \left \{ \omega \in \Omega:\left( f_{n}\left( \omega \right)
,...,f_{n+k}\left( \omega \right) \right) \in B\right \} \right) =\nu \left(
\left \{ \omega \in \Omega:\left( f_{n+1}\left( \omega \right)
,...,f_{n+k+1}\left( \omega \right) \right) \in B\right \} \right) .
\label{eq:sta}
\end{equation}
\end{definition}

This notion generalizes the usual notion of stationary stochastic process by
allowing for the nonadditivity of the underlying probability measure. Recall
that $\left( 
\mathbb{R}
^{%
\mathbb{N}
},\sigma \left( \mathcal{C}\right) \right) $ denotes the space of sequences
endowed with the $\sigma $-algebra generated by the algebra of cylinders. We
denote a generic element of $%
\mathbb{R}
^{%
\mathbb{N}
}$ by $x$. We also consider the shift transformation $\tau :%
\mathbb{R}
^{%
\mathbb{N}
}\rightarrow 
\mathbb{R}
^{%
\mathbb{N}
}$ defined by%
\begin{equation*}
\tau \left( x\right) =\left( x_{2},x_{3},x_{4},......\right) \qquad \forall
x\in 
\mathbb{R}
^{%
\mathbb{N}
}.
\end{equation*}%
The sequence $\left\{ f_{n}\right\} _{n\in 
\mathbb{N}
}$ induces a natural (measurable) map between $\left( \Omega ,\mathcal{F}%
\right) $ and $\left( 
\mathbb{R}
^{%
\mathbb{N}
},\sigma \left( \mathcal{C}\right) \right) $, defined by%
\begin{equation*}
\omega \mapsto \mathbf{f}\left( \omega \right) =\left( f_{1}\left( \omega
\right) ,...,f_{k}\left( \omega \right) ,...\right) \qquad \forall \omega
\in \Omega .
\end{equation*}%
Define $\nu _{\mathbf{f}}:\sigma \left( \mathcal{C}\right) \rightarrow \left[
0,1\right] $ by%
\begin{equation*}
\nu _{\mathbf{f}}\left( C\right) =\nu \left( \mathbf{f}^{-1}\left( C\right)
\right) \qquad \forall C\in \sigma \left( \mathcal{C}\right) .
\end{equation*}

\begin{definition}
Given a capacity $\nu$, we say that $\mathbf{f}=\left \{ f_{n}\right \}
_{n\in%
\mathbb{N}
}$ is ergodic if and only if $\nu_{\mathbf{f}}$ is ergodic with respect to
the shift transformation.
\end{definition}

\begin{lemma}
\label{lem:law}If $\nu $ is a convex capacity continuous at $\Omega $ and $%
\mathbf{f}$ is stationary, then $\nu _{\mathbf{f}}$ is a convex capacity
continuous at $%
\mathbb{R}
^{%
\mathbb{N}
}$ which is shift invariant. Moreover, $\mathbf{f}$ is ergodic if $\nu
\left( \mathcal{T}\right) =\left \{ 0,1\right \} $.
\end{lemma}

This observation is a first step to deduce the Strong Law of Large Numbers
as a corollary of Theorem \ref{thm:w-inv-erg} applied to $\nu _{\mathbf{f}}$%
. In a nutshell, the assumption of stationarity yields that the limit%
\begin{equation*}
\lim_{n}\frac{1}{n}\sum_{k=1}^{n}f_{k}
\end{equation*}%
exists $\nu $-a.s. In order to obtain also a characterization of the limit
in terms of the (Choquet) expected value, we further need $\nu _{\mathbf{f}}$
to be ergodic.

\begin{theorem}
\label{thm:SLLN}Let $\nu $ be a convex capacity continuous at $\Omega $. If $%
\mathbf{f}=\left \{ f_{n}\right \} _{n\in 
\mathbb{N}
}$ is stationary and ergodic, then%
\begin{equation*}
\nu \left( \left \{ \omega \in \Omega :\int_{\Omega }f_{1}d\nu \leq \lim_{n}%
\frac{1}{n}\sum_{k=1}^{n}f_{k}\left( \omega \right) \leq \int_{\Omega }f_{1}d%
\bar{\nu}\right \} \right) =1.
\end{equation*}
\end{theorem}

We close by observing that there are few but important differences with the
nonadditive Strong Law of Large Numbers of Marinacci \cite{Mar}\ and
Maccheroni and Marinacci \cite{MMa}. In terms of hypotheses, we weaken the
assumption of total monotonicity of $\nu $ to convexity, while we replace
the i.i.d hypothesis of \cite{Mar} with stationarity and ergodicity.
Finally, compared to the main result of \cite{MMa}, we need to assume the
continuity of $\nu $. In turn, we obtain that empirical averages exist $\nu $%
-a.s., a property that was not present in previous works. The bounds for
these empirical averages are in terms of the lower and the upper Choquet
integrals of the random variable $f_{1}$, as in \cite{Mar}\ and \cite{MMa}.

\appendix

\section{Dynkin Spaces and Nonadditive Probabilities}

Consider a standard measurable space $\left( \Omega ,\mathcal{F}\right) $
and a transformation $\tau :\Omega \rightarrow \Omega $ which is $\mathcal{F}%
\backslash \mathcal{F}$ measurable. Recall that we denote by $\mathcal{I}$
the set of all invariant probability measures. If $\mathcal{I}$ is a
nonempty set, then the triple $\left( \Omega ,\mathcal{F},\mathcal{I}\right) 
$ forms a Dynkin space.

\begin{definition}[Dynkin, 1978]
\label{def:Dyn}Let $\mathcal{P}\ $be a nonempty subset of $%
\Delta^{\sigma}\left( \Omega,\mathcal{F}\right) $ where $\left( \Omega,%
\mathcal{F}\right) $\ is a separable measurable space. The triple $\left(
\Omega,\mathcal{F},\mathcal{P}\right) $ is a Dynkin space if and only if
there exist a sub-$\sigma$-algebra $\mathcal{G}\subseteq \mathcal{F}$, a set 
$W\in \mathcal{F}$, and a function%
\begin{equation*}
\begin{array}{cccc}
p: & \mathcal{F}\times \Omega & \rightarrow & \left[ 0,1\right] \\ 
& \left( A,\omega \right) & \mapsto & p\left( A,\omega \right)%
\end{array}%
\end{equation*}
such that:

\begin{enumerate}
\item[(a)] for each $P\in \mathcal{P}$ and $A\in \mathcal{F}$, $p\left(
A,\cdot \right) :\Omega \rightarrow \left[ 0,1\right] $ is a version of the
conditional probability of $A$ given $\mathcal{G}$;

\item[(b)] for each $\omega \in \Omega$, $p\left( \cdot,\omega \right) :%
\mathcal{F}\rightarrow \left[ 0,1\right] $ is a probability measure;

\item[(c)] $P\left( W\right) =1$ for all $P\in \mathcal{P}$ and $p\left(
\cdot,\omega \right) \in \mathcal{P}$ for all $\omega \in W$.
\end{enumerate}
\end{definition}

It is not hard to check that, given $f\in B\left( \Omega ,\mathcal{F}\right) 
$, the function $\hat{f}:\Omega \rightarrow 
\mathbb{R}
$, defined by%
\begin{equation}
\hat{f}\left( \omega \right) =\int_{\Omega }fdp\left( \cdot ,\omega \right)
\qquad \forall \omega \in \Omega ,  \label{eq:cond}
\end{equation}%
is a version of the conditional expected value of $f$ given $\mathcal{G}$
for all $P\in \mathcal{P}$, in particular, $\hat{f}\in B\left( \Omega ,%
\mathcal{G}\right) $ (see also \cite[Remark 13]{CMMM6}). When $\left( \Omega
,\mathcal{F}\right) $ is a standard measurable space, if $\left( \Omega ,%
\mathcal{F},\mathcal{P}\right) =\left( \Omega ,\mathcal{F},\mathcal{I}%
\right) $, then $\mathcal{G}$ is the set of invariant events. In particular,
we can consider $W=\Omega $ (see Gray \cite[Theorem 8.3]{Gra}). We conclude
with an ancillary lemma.

\begin{lemma}
\label{lem:erg}Let $\left( \Omega ,\mathcal{F}\right) $ be a measurable
space and $\mathcal{G}$ a sub-$\sigma $-algebra of $\mathcal{F}$. If $\nu $
is a lower probability such that $\nu \left( \mathcal{G}\right) =\left \{
0,1\right \} $ and $g\in B\left( \Omega ,\mathcal{G}\right) $, then%
\begin{equation*}
\nu \left( \left \{ \omega \in \Omega :\int_{\Omega }gd\nu \leq g\left(
\omega \right) \leq \int_{\Omega }gd\bar{\nu}\right \} \right) =1.
\end{equation*}
\end{lemma}

\noindent \textbf{Proof. }We proceed by assuming that $g\geq 0$. Since $\nu $
is a capacity such that $\nu \left( \mathcal{G}\right) =\left \{
0,1\right
\} $ and $0\leq g\leq \lambda $ for some $\lambda \in 
\mathbb{R}
$, it follows that the sets%
\begin{align*}
I& =\left \{ t\in \left[ 0,\infty \right) :\nu \left( \left \{ \omega \in
\Omega :g\left( \omega \right) \geq t\right \} \right) =1\right \} \\
& \text{and} \\
J& =\left \{ t\in \left( -\infty ,0\right] :\nu \left( \left \{ \omega \in
\Omega :-g\left( \omega \right) \geq t\right \} \right) =1\right \}
\end{align*}%
are well defined nonempty intervals. $I$ is bounded from above and such that 
$0\in I$. $J$ is unbounded from below and such that $-\lambda \in J$. Since $%
\nu $ is a lower probability, $\nu $ is continuous. We can conclude that $%
t^{\star }=\sup I\in I$ and $t_{\star }=\sup J\in J$. Since $\nu \left( 
\mathcal{G}\right) =\left \{ 0,1\right \} $, this implies that%
\begin{align*}
\int_{\Omega }gd\nu & =\int_{0}^{\infty }\nu \left( \left \{ \omega \in
\Omega :g\left( \omega \right) \geq t\right \} \right) dt=\int_{0}^{\sup
I}dt=t^{\star } \\
& and \\
\int_{\Omega }-gd\nu & =\int_{-\infty }^{0}\left[ \nu \left( \left \{ \omega
\in \Omega :-g\left( \omega \right) \geq t\right \} \right) -\nu \left(
\Omega \right) \right] dt=\int_{\sup J}^{0}\left( -1\right) dt=t_{\star }.
\end{align*}%
It follows that $t^{\star }=\int_{\Omega }gd\nu $ and $t_{\star
}=\int_{\Omega }-gd\nu $. Since $t^{\star }\in I$ and $t_{\star }\in J$, we
also have that%
\begin{equation*}
\nu \left( \left \{ \omega \in \Omega :g\left( \omega \right) \geq t^{\star
}\right \} \right) =1=\nu \left( \left \{ \omega \in \Omega :g\left( \omega
\right) \leq -t_{\star }\right \} \right) .
\end{equation*}%
Since $\nu $ is a lower probability, this implies that%
\begin{equation}
\nu \left( \left \{ \omega \in \Omega :\int_{\Omega }gd\nu \leq g\left(
\omega \right) \leq \int_{\Omega }gd\bar{\nu}\right \} \right) =\nu \left(
\left \{ \omega \in \Omega :t^{\star }\leq g\left( \omega \right) \leq
-t_{\star }\right \} \right) =1.  \label{eq:pos-san}
\end{equation}%
We next remove the hypothesis that $g\geq 0$. Since $g\in B\left( \Omega ,%
\mathcal{G}\right) $, it follows that there exists $c\in 
\mathbb{R}
$\ such that $g+c1_{\Omega }\geq 0$. By (\ref{eq:pos-san}) and since the
Choquet integral is constant additive, it follows that%
\begin{align*}
1& =\nu \left( \left \{ \omega \in \Omega :\int_{\Omega }\left( g+c1_{\Omega
}\right) d\nu \leq g\left( \omega \right) +c\leq \int_{\Omega }\left(
g+c1_{\Omega }\right) d\bar{\nu}\right \} \right) \\
& =\nu \left( \left \{ \omega \in \Omega :\int_{\Omega }gd\nu +c\leq g\left(
\omega \right) +c\leq \int_{\Omega }gd\bar{\nu}+c\right \} \right) \\
& =\nu \left( \left \{ \omega \in \Omega :\int_{\Omega }gd\nu \leq g\left(
\omega \right) \leq \int_{\Omega }gd\bar{\nu}\right \} \right) ,
\end{align*}%
proving the statement.\hfill $\blacksquare $

\section{Proofs}

\noindent \textbf{Proof of Proposition \ref{pro:inv-low}. }Recall that if $%
\nu $ is a lower probability, we have that%
\begin{equation}
\nu \leq P\leq \bar{\nu}\qquad \forall P\in \mathrm{core}\left( \nu \right)
\subseteq \Delta ^{\sigma }\left( \Omega ,\mathcal{F}\right) .
\label{eq:sand-low}
\end{equation}

\smallskip

1. Pick $A\in \mathcal{F}$. Since $\nu $ is strongly invariant and $\nu \leq 
\bar{\nu}$, we have $\bar{\nu}\left( \tau ^{-1}\left( A\right) \backslash
A\right) =\nu \left( A\backslash \tau ^{-1}\left( A\right) \right) \leq \bar{%
\nu}\left( A\backslash \tau ^{-1}\left( A\right) \right) =\nu \left( \tau
^{-1}\left( A\right) \backslash A\right) \leq \bar{\nu}\left( \tau
^{-1}\left( A\right) \backslash A\right) $. It follows that\ $\nu \left(
A\backslash \tau ^{-1}\left( A\right) \right) =\bar{\nu}\left( A\backslash
\tau ^{-1}\left( A\right) \right) =\bar{\nu}\left( \tau ^{-1}\left( A\right)
\backslash A\right) =\nu \left( \tau ^{-1}\left( A\right) \backslash
A\right) =k$. By (\ref{eq:sand-low}), we can conclude that $P\left(
A\backslash \tau ^{-1}\left( A\right) \right) =k=P\left( \tau ^{-1}\left(
A\right) \backslash A\right) $ for all $P\in \mathrm{core}\left( \nu \right) 
$. This implies that $P\left( A\right) =P\left( A\backslash \tau ^{-1}\left(
A\right) \right) +P\left( A\cap \tau ^{-1}\left( A\right) \right) =P\left(
\tau ^{-1}\left( A\right) \backslash A\right) +P\left( A\cap \tau
^{-1}\left( A\right) \right) =P\left( \tau ^{-1}\left( A\right) \right) $
for all $P\in \mathrm{core}\left( \nu \right) $, proving the statement.

\smallskip

2. By assumption, there exists a convex capacity $\rho :\mathcal{A}_{%
\mathcal{S}\left( \mathcal{I}\right) }\rightarrow \left[ 0,1\right] $ such
that%
\begin{equation}
\nu \left( A\right) =\int_{\mathcal{S}\left( \mathcal{I}\right) }P\left(
A\right) d\rho \left( P\right) =\min_{\pi \in \mathrm{core}\left( \rho
\right) }\int_{\mathcal{S}\left( \mathcal{I}\right) }P\left( A\right) d\pi
\left( P\right) \qquad \forall A\in \mathcal{F}.  \label{eq:rob-pro}
\end{equation}%
Define $\mathcal{M}=\left \{ \nu _{\pi }:\pi \in \mathrm{core}\left( \rho
\right) \right \} $. By \cite[Lemma 24]{CMMM6} and (\ref{eq:rob-pro}) and
since $\nu $ is continuous at $\Omega $, we have that $\rho $ is continuous
at $\mathcal{S}\left( \mathcal{I}\right) $, thus, each $\pi $ in $\mathrm{%
core}\left( \rho \right) $ is a probability measure and $\mathcal{M}$ is a
compact subset of $\Delta ^{\sigma }\left( \Omega ,\mathcal{F}\right) $.
Moreover, we also have that $\mathcal{M}\subseteq \mathcal{I}$. We can
conclude that $\nu \left( A\right) =\min_{\pi \in \mathrm{core}\left( \rho
\right) }\int_{\mathcal{S}\left( \mathcal{I}\right) }P\left( A\right) d\pi
\left( P\right) =\min_{P\in \mathcal{M}}P\left( A\right) $ for all $A\in 
\mathcal{F}$, proving the statement.

$\smallskip $

3. Fix $\mathcal{M}\in \mathcal{A}_{\mathcal{S}\left( \mathcal{I}\right) }$.
Consider $\rho :\mathcal{A}_{\mathcal{S}\left( \mathcal{I}\right)
}\rightarrow \left[ 0,1\right] $ defined by%
\begin{equation*}
\rho \left( F\right) =\left \{ 
\begin{array}{cc}
1 & F\supseteq \mathcal{M} \\ 
0 & otherwise%
\end{array}%
\right. \qquad \forall F\in \mathcal{A}_{\mathcal{S}\left( \mathcal{I}%
\right) }.
\end{equation*}%
It is immediate to check that $\rho $ is a convex capacity. By \cite[Example
4.4]{MM} and since $\mathcal{M}\in \mathcal{A}_{\mathcal{S}\left( \mathcal{I}%
\right) }$, we have that $\nu \left( A\right) =\min_{P\in \mathcal{M}%
}P\left( A\right) =\int_{\mathcal{S}\left( \mathcal{I}\right) }P\left(
A\right) d\rho \left( P\right) $ for all $A\in \mathcal{F}$. Since $\mathcal{%
M}\subseteq \mathcal{S}\left( \mathcal{I}\right) $, observe that $P\left(
A\right) \in \left \{ 0,1\right \} $ for all $P\in \mathcal{M}$ and for all $%
A\in \mathcal{G}$. It follows that $\nu \left( \mathcal{G}\right) =\left \{
0,1\right \} $.

\smallskip

4. Since $\nu $ is a functionally invariant lower probability, we have that $%
\mathcal{M}\subseteq \mathcal{I}$ and $\nu \left( A\right) =\min_{P\in 
\mathcal{M}}P\left( A\right) =\min_{P\in \mathcal{M}}P\left( \tau
^{-1}\left( A\right) \right) =\nu \left( \tau ^{-1}\left( A\right) \right) $
for all $A\in \mathcal{F}$, proving that $\nu $ is invariant.\hfill $%
\blacksquare $

\smallskip

\noindent \textbf{Proof of Theorem \ref{thm:inv-low-con}. }Recall that if $%
\nu $ is convex and continuous at $\Omega$, then it is a lower probability.

\smallskip

(i) implies (ii). It follows by point 1 of Proposition \ref{pro:inv-low}.

\smallskip

(ii) implies (iii). We just need to show that $\nu $ is robustly invariant.
Define $I:B\left( \Omega ,\mathcal{F}\right) \rightarrow 
\mathbb{R}
$ by%
\begin{equation*}
I\left( f\right) =\int_{\Omega }fd\nu \qquad \forall f\in B\left( \Omega ,%
\mathcal{F}\right) .
\end{equation*}%
By Schmeidler \cite{Sch1} (see also \cite{MM}), $I$ is comonotonic additive
and supermodular. Since $\nu $ is convex, we have that $I\left( f\right)
=\min_{P\in \mathrm{core}\left( \nu \right) }\int_{\Omega }fdP$ for all $%
f\in B\left( \Omega ,\mathcal{F}\right) $. Since $\mathrm{core}\left( \nu
\right) \subseteq \mathcal{I}$, this implies that if $\int_{\Omega
}f_{1}dP\geq \int_{\Omega }f_{2}dP$ for all $P\in \mathcal{I}$, then $%
I\left( f_{1}\right) \geq I\left( f_{2}\right) $. In particular, $I\left(
f\right) =I\left( \hat{f}\right) $ for all $f\in B\left( \Omega ,\mathcal{F}%
\right) $. It is also immediate to see that $I\left( k1_{\Omega }\right) =k$
for all $k\in 
\mathbb{R}
$. It follows that $I$ restricted to $B\left( \Omega ,\mathcal{G}\right) $
is normalized, comonotonic additive, supermodular, and such that $%
\int_{\Omega }f_{1}dP\geq \int_{\Omega }f_{2}dP$ for all $P\in \mathcal{I}$
implies $I\left( f_{1}\right) \geq I\left( f_{2}\right) $. By \cite[Lemma 24
and Proposition 25]{CMMM6} and since $\left( \Omega ,\mathcal{F},\mathcal{I}%
\right) $ is a Dynkin space, it follows that there exists $\breve{I}:B\left( 
\mathcal{S}\left( \mathcal{I}\right) ,\mathcal{A}_{\mathcal{S}\left( 
\mathcal{I}\right) }\right) \rightarrow 
\mathbb{R}
$ such that $\breve{I}$ is normalized, monotone, comonotonic additive,
supermodular, and such that $I\left( f\right) =\breve{I}\left( \left\langle
f,\cdot \right\rangle \right) $ for all $f\in B\left( \Omega ,\mathcal{G}%
\right) $. By \cite{Sch1} (see also \cite{MM}), it follows that there exists
a convex capacity $\rho :\mathcal{A}_{\mathcal{S}\left( \mathcal{I}\right)
}\rightarrow \left[ 0,1\right] $ such that%
\begin{equation}
I\left( f\right) =\int_{\mathcal{S}\left( \mathcal{I}\right) }\left(
\int_{\Omega }fdP\right) d\rho \left( P\right) \qquad \forall f\in B\left(
\Omega ,\mathcal{G}\right) .  \label{eq:pomp}
\end{equation}%
Since $I\left( f\right) =I\left( \hat{f}\right) $ for all $f\in B\left(
\Omega ,\mathcal{F}\right) $, it follows that (\ref{eq:pomp}) holds for all $%
f\in B\left( \Omega ,\mathcal{F}\right) $. In particular, by picking $%
f=1_{A} $ with $A\in \mathcal{F}$, this shows that $\nu $ is robustly
invariant.

$\smallskip$

(iii) implies (iv). It is trivial.

\smallskip

(iv) implies (i). Since $\nu $ is convex and $\mathrm{core}\left( \nu
\right) \subseteq \mathcal{I}$, it follows that 
\begin{eqnarray*}
\nu \left( A\backslash \tau ^{-1}\left( A\right) \right) +\nu \left( A\cup
\left( \tau ^{-1}\left( A\right) \right) ^{c}\right) &=&\int_{\Omega }\left(
1_{\Omega }+1_{A}-1_{\tau ^{-1}\left( A\right) }\right) d\nu \\
&=&\min_{P\in \mathrm{core}\left( \nu \right) }\int_{\Omega }\left(
1_{\Omega }+1_{A}-1_{\tau ^{-1}\left( A\right) }\right) dP=1.
\end{eqnarray*}%
Thus,\ we have that%
\begin{equation*}
\nu \left( A\backslash \tau ^{-1}\left( A\right) \right) =1-\nu \left( A\cup
\left( \tau ^{-1}\left( A\right) \right) ^{c}\right) =1-\nu \left( \left(
\tau ^{-1}\left( A\right) \backslash A\right) ^{c}\right) =\bar{\nu}\left(
\tau ^{-1}\left( A\right) \backslash A\right) .
\end{equation*}%
An analogous argument yields that $\nu \left( \tau ^{-1}\left( A\right)
\backslash A\right) =\bar{\nu}\left( A\backslash \tau ^{-1}\left( A\right)
\right) $, proving the statement.\hfill $\blacksquare $

\smallskip

Before proving Theorem \ref{thm:w-inv-erg}, we provide an ancillary key
result.

\begin{theorem}
\label{thm:erg-PI}Let $\left( \Omega ,\mathcal{F}\right) $ be a measurable
space, $\nu $ a lower probability, and assume that the family $\mathcal{I}$
of invariant probability measures is not empty. The following statements are
equivalent:

\begin{enumerate}
\item[(i)] There exists $\breve{P}\in \mathcal{I}$ such that for each $E\in%
\mathcal{F}$%
\begin{equation*}
\breve{P}\left( E\right) =1\implies \lim_{k}\nu \left( \tau^{-k}\left(
E\right) \right) =1;
\end{equation*}

\item[(ii)] There exists $\breve{P}\in \mathcal{I}$ such that for each $E\in%
\mathcal{G}$%
\begin{equation*}
\breve{P}\left( E\right) =1\implies \nu \left( E\right) =1;
\end{equation*}

\item[(iii)] For each $E\in \mathcal{G}$%
\begin{equation*}
P\left( E\right) =1\qquad \forall P\in \mathcal{I}\implies \nu \left(
E\right) =1;
\end{equation*}

\item[(iv)] For each $f\in B\left( \Omega ,\mathcal{F}\right) $ there exists 
$f^{\star }\in B\left( \Omega ,\mathcal{G}\right) $ such that%
\begin{equation*}
\lim_{n}\frac{1}{n}\dsum \limits_{k=1}^{n}f\left( \tau ^{k-1}\left( \omega
\right) \right) =f^{\star }\left( \omega \right) \qquad \nu -a.s.;
\end{equation*}

\item[(v)] \textrm{core}$\left( \nu \right) \subseteq \mathcal{PI}$.
\end{enumerate}
\end{theorem}

\noindent \textbf{Proof. }(i) implies (ii). If $E\in \mathcal{G}$, then $%
\tau^{-k}\left( E\right) =E$ for all $k\in%
\mathbb{N}
$, yielding the statement.

\smallskip

(ii) implies (iii). It is trivial.

\smallskip

(iii) implies (iv). Consider $f\in B\left( \Omega,\mathcal{F}\right) $.
Define $f^{\star}:\Omega \rightarrow%
\mathbb{R}
$ by%
\begin{equation*}
f^{\star}\left( \omega \right) =\limsup_{n}\frac{1}{n}\dsum
\limits_{k=1}^{n}f\left( \tau^{k-1}\left( \omega \right) \right) \qquad
\forall \omega \in \Omega.
\end{equation*}
Define $f_{\star}:\Omega \rightarrow%
\mathbb{R}
$ by considering the $\liminf$. Since $f\in B\left( \Omega,\mathcal{F}%
\right) $, it can be shown that $f^{\star},f_{\star}\in B\left( \Omega,%
\mathcal{G}\right) $. Consider the event 
\begin{align*}
E & =\left \{ \omega \in \Omega:\lim_{n}\frac{1}{n}\dsum
\limits_{k=1}^{n}f\left( \tau^{k-1}\left( \omega \right) \right) \text{
exists}\right \} =\left \{ \omega \in \Omega:f^{\star}\left( \omega \right)
=f_{\star}\left( \omega \right) \right \} \\
& =\left \{ \omega \in \Omega:f^{\star}\left( \omega \right) =\lim_{n}\frac {%
1}{n}\dsum \limits_{k=1}^{n}f\left( \tau^{k-1}\left( \omega \right) \right)
=f_{\star}\left( \omega \right) \right \} .
\end{align*}
By Birkhoff's Ergodic Theorem (see \cite[Theorem 24.1]{Bil}), we have that $%
P\left( E\right) =1$ for all $P\in \mathcal{I}$. By assumption, this yields
that $\nu \left( E\right) =1$. Since $f$ was chosen to be generic, the
statement follows.

\smallskip

(iv) implies (v). Recall that for each $P\in \mathrm{core}\left( \nu \right) 
$, $P\left( A\right) \geq \nu \left( A\right) $ for all $A\in \mathcal{F}$.%
\textrm{\ }By assumption, we can conclude that for each $P\in \mathrm{core}%
\left( \nu \right) $, for each $f\in B\left( \Omega ,\mathcal{F}\right) $
there exists $f^{\star }\in B\left( \Omega ,\mathcal{G}\right) $ such that%
\begin{equation*}
\lim_{n}\frac{1}{n}\dsum\limits_{k=1}^{n}f\left( \tau ^{k-1}\left( \omega
\right) \right) =f^{\star }\left( \omega \right) \qquad P-a.s.
\end{equation*}%
By \cite[p. 964]{GK} (see also \cite[Exercises 31 and 32, pp. 723--724]{DS}%
), it follows that $P\in \mathcal{PI}$.

\smallskip

(v) implies (i). Since $\nu $ is a lower probability, it is continuous at $%
\Omega $ and exact. By \cite[Theorem 4.2]{MM}, it follows that there exists
a measure $P\in \mathrm{core}\left( \nu \right) $ such that for each $A\in 
\mathcal{F}$, for each $\varepsilon >0$, there exists $\delta >0$ such that%
\begin{equation}
P\left( A\right) <\delta \Longrightarrow Q\left( A\right) <\varepsilon
\qquad \forall Q\in \mathrm{core}\left( \nu \right) .  \label{eq:dom-PI}
\end{equation}%
It is immediate to show that $P$ is such that for each $A\in \mathcal{F}$%
\begin{equation}
P\left( A\right) =0\Longrightarrow Q\left( A\right) =0\qquad \forall Q\in 
\mathrm{core}\left( \nu \right) .  \label{eq:con}
\end{equation}%
Since $P\in \mathrm{core}\left( \nu \right) \subseteq \mathcal{PI}$, we have
that there exists $\breve{P}\in \mathcal{I}$ such that $\breve{P}\left(
E\right) =P\left( E\right) $ for all $E\in \mathcal{G}$. Consider $E\in 
\mathcal{F}$. Assume that $\breve{P}\left( E\right) =1$. It follows that $%
\breve{P}\left( E^{c}\right) =0$. At the same time, define $F_{n}=\cup
_{k=n}^{\infty }\tau ^{-k}\left( E^{c}\right) $. Note that $F_{n}\downarrow
F\in \mathcal{G}$. Since $\breve{P}\in \mathcal{I}$, it follows that $\breve{%
P}\left( F\right) =\lim_{n}\breve{P}\left( F_{n}\right) \leq \breve{P}\left(
F_{1}\right) \leq \sum_{k=1}^{\infty }\breve{P}\left( \tau ^{-k}\left(
E^{c}\right) \right) =0$. It follows that $\breve{P}\left( F\right) =0$,
that is, $P\left( F\right) =0$. By (\ref{eq:con}), we have that $Q\left(
F\right) =0$ for all $Q\in \mathrm{core}\left( \nu \right) $, that is, $\bar{%
\nu}\left( F\right) =0$. Since $\nu $ is a lower probability, $\bar{\nu}$
satisfies the Fatou's property, that is, $0\leq \limsup_{k}\bar{\nu}\left(
A_{k}\right) \leq \bar{\nu}\left( \limsup_{k}A_{k}\right) $ for each\
sequence $\left \{ A_{k}\right \} _{k\in 
\mathbb{N}
}\subseteq \mathcal{F}$. This implies that $0\leq \liminf_{k}\bar{\nu}\left(
\tau ^{-k}\left( E^{c}\right) \right) \leq \limsup_{k}\bar{\nu}\left( \tau
^{-k}\left( E^{c}\right) \right) \leq \bar{\nu}\left( \limsup_{k}\tau
^{-k}\left( E^{c}\right) \right) =\bar{\nu}\left( F\right) =0$. We can
conclude that $\lim_{k}\nu \left( \tau ^{-k}\left( E\right) \right) =\lim_{k}%
\left[ 1-\bar{\nu}\left( \tau ^{-k}\left( E^{c}\right) \right) \right] =1$,
proving the statement.\hfill $\blacksquare $

\smallskip

The proof of Theorem \ref{thm:w-inv-erg} uses some of the techniques common
in Ergodic Theory (see, e.g., \cite[Theorem 7]{Dow1}). Also, note that,
given a capacity $\nu$, we have that%
\begin{equation*}
\mathrm{core}\left( \nu \right) =\left \{ P\in \Delta \left( \Omega ,%
\mathcal{F}\right) :\bar{\nu}\geq P\geq \nu \right \} =\left \{ P\in \Delta
\left( \Omega,\mathcal{F}\right) :\bar{\nu}\geq P\right \} .
\end{equation*}

\smallskip

\noindent \textbf{Proof of Theorem \ref{thm:w-inv-erg}. }We first prove
that, given the assumptions, $\emptyset \not =$\textrm{core}$\left( \nu
\right) \subseteq \mathcal{PI}$. In particular, this shows that $\mathcal{I}%
\not =\emptyset$.

\smallskip

\textit{Claim: Let }$\nu$\textit{\ be a lower probability. If }$\nu$\textit{%
\ is invariant, then }$\mathrm{core}\left( \nu \right) \subseteq \mathcal{PI}
$\textit{. In particular, }$\mathcal{I}\not =\emptyset$\textit{.}

\smallskip

\textit{Proof of the Claim. }Since $\nu $ is invariant, $\bar{\nu}$ is
invariant. Since $\nu $ is a lower probability, $\nu $ is continuous at $%
\Omega $ and, in particular, $\emptyset \not=\mathrm{core}\left( \nu \right)
\subseteq \Delta ^{\sigma }\left( \Omega ,\mathcal{F}\right) $. Fix a
Banach-Mazur limit (see \cite[pag. 550]{AB}) $\phi :l^{\infty }\rightarrow 
\mathbb{R}
$, that is, a functional from $l^{\infty }$ to $%
\mathbb{R}
$\ such that:

\begin{enumerate}
\item $\phi$ is linear;

\item $\phi$ is positive;

\item $\phi \left( x_{1},x_{2},...\right) =\phi \left( x_{2},x_{3}...\right) 
$ for all $x\in l^{\infty}$;

\item $\phi \left( x_{1},x_{2},...\right) =\lim_{n}x_{n}$ for all $x\in c$.
\end{enumerate}

Observe that $\nu \left( A\right) \leq P\left( A\right) \leq \bar{\nu}\left(
A\right) $ for all $P\in \mathrm{core}\left( \nu \right) $ and all $A\in 
\mathcal{F}$. Fix $P\in \mathrm{core}\left( \nu \right) $, define $P_{n}:%
\mathcal{F}\rightarrow \left[ 0,1\right] $ by%
\begin{equation*}
P_{n}\left( A\right) =\frac{1}{n}\sum_{k=0}^{n-1}P\left( \tau ^{-k}\left(
A\right) \right) \qquad \forall A\in \mathcal{F}.
\end{equation*}%
Note that $P\left( \tau ^{-k}\left( A\right) \right) \leq \bar{\nu}\left(
\tau ^{-k}\left( A\right) \right) =\bar{\nu}\left( A\right) $ for all $A\in 
\mathcal{F}$ and for all $k\in 
\mathbb{N}
_{0}$. Since $\mathrm{core}\left( \nu \right) $ is convex, this implies that 
$\left \{ P_{n}\right \} _{n\in 
\mathbb{N}
}\subseteq \mathrm{core}\left( \nu \right) $. For each $A\in \mathcal{F}$,
define\ $x_{A}=\left( P_{1}\left( A\right) ,P_{2}\left( A\right)
,P_{3}\left( A\right) ,...\right) $. Note that $0\leq x_{A}\leq 1_{%
\mathbb{N}
}$, thus, $x_{A}\in l^{\infty }$ for all $A\in \mathcal{F}$. Define $\hat{P}:%
\mathcal{F}\rightarrow \left[ 0,1\right] $ by $\hat{P}\left( A\right) =\phi
\left( x_{A}\right) $ for all $A\in \mathcal{F}$. Since $\phi $ is positive,
note that $\hat{P}$ is a well defined positive set function. Next, consider $%
A,B\in \mathcal{F}$ such that $A\cap B=\emptyset $. Since $\left \{
P_{n}\right \} _{n\in 
\mathbb{N}
}\subseteq \Delta \left( \Omega ,\mathcal{F}\right) $, it follows that $%
P_{n}\left( A\cup B\right) =P_{n}\left( A\right) +P_{n}\left( B\right) $ for
all $n\in 
\mathbb{N}
$. Since $\phi $ is linear, this implies that $\hat{P}\left( A\cup B\right)
=\phi \left( x_{A\cup B}\right) =\phi \left( x_{A}+x_{B}\right) =\phi \left(
x_{A}\right) +\phi \left( x_{B}\right) =\hat{P}\left( A\right) +\hat{P}%
\left( B\right) $, proving that $\hat{P}$ is additive. Next, consider $A\in 
\mathcal{G}$. Since $\tau ^{-k}\left( A\right) =A$ for all $k\in 
\mathbb{N}
_{0}$, it follows that $P_{n}\left( A\right) =P\left( A\right) $ for all $%
n\in 
\mathbb{N}
$. Since $\phi $ maps convergent sequences into their limit, we have that $%
\hat{P}\left( A\right) =\phi \left( x_{A}\right) =P\left( A\right) $. In
particular, this implies that $\hat{P}\left( \Omega \right) =1$ and $\hat{P}%
\left( \emptyset \right) =0$. Up to now, we have proved that $\hat{P}\in
\Delta \left( \Omega ,\mathcal{F}\right) $ and $\hat{P}\left( A\right)
=P\left( A\right) $ for all $A\in \mathcal{G}$. Since $\left \{
P_{n}\right
\} _{n\in 
\mathbb{N}
}\subseteq \mathrm{core}\left( \nu \right) $, we have that $x_{A}\leq \bar{%
\nu}\left( A\right) 1_{%
\mathbb{N}
}$. Since $\phi $ is linear and positive, it follows that $\hat{P}\left(
A\right) =\phi \left( x_{A}\right) \leq \phi \left( \bar{\nu}\left( A\right)
1_{%
\mathbb{N}
}\right) =\bar{\nu}\left( A\right) $ for all $A\in \mathcal{F}$, that is, $%
\hat{P}\in \mathrm{core}\left( \nu \right) $. Since $\mathrm{core}\left( \nu
\right) \subseteq \Delta ^{\sigma }\left( \Omega ,\mathcal{F}\right) $, we
can conclude that $\hat{P}\in \Delta ^{\sigma }\left( \Omega ,\mathcal{F}%
\right) $. We next show that $\hat{P}$ is invariant. Note that for each $%
A\in \mathcal{F}$ and for each $n\in 
\mathbb{N}
$%
\begin{align*}
P_{n}\left( \tau ^{-1}\left( A\right) \right) & =\frac{1}{n}%
\sum_{k=0}^{n-1}P\left( \tau ^{-k-1}\left( A\right) \right) =\frac{n+1}{n}%
\cdot \frac{1}{n+1}\sum_{k=0}^{n}P\left( \tau ^{-k}\left( A\right) \right) -%
\frac{1}{n}P\left( A\right) \\
& =\frac{n+1}{n}P_{n+1}\left( A\right) -\frac{1}{n}P\left( A\right) .
\end{align*}%
Define $y=\left( P_{2}\left( A\right) ,P_{3}\left( A\right) ,...\right) $.
Define $z=x_{\tau ^{-1}\left( A\right) }-y\in l^{\infty }$. Note that $%
\left
\vert z_{n}\right \vert =\left \vert P_{n}\left( \tau ^{-1}\left(
A\right) \right) -P_{n+1}\left( A\right) \right \vert \leq \frac{1}{n}%
\left
\vert P_{n+1}\left( A\right) -P\left( A\right) \right \vert \leq 
\frac{2}{n}$ for all $n\in 
\mathbb{N}
$. It follows that $\lim_{n}z_{n}=0$. Since $\phi $ satisfies properties 3,
1, and 4, we have that $\left \vert \hat{P}\left( \tau ^{-1}\left( A\right)
\right) -\hat{P}\left( A\right) \right \vert =\left \vert \phi \left(
x_{\tau ^{-1}\left( A\right) }\right) -\phi \left( x_{A}\right) \right \vert
=\left \vert \phi \left( x_{\tau ^{-1}\left( A\right) }\right) -\phi \left(
y\right) \right \vert =\left \vert \phi \left( z\right) \right \vert =0$,
proving that $\hat{P}$ is invariant. Given the previous part of the proof, $%
\hat{P}\in \mathcal{I}$ and $P\in \mathcal{PI}$. Since $P$ was arbitrarily
chosen in $\mathrm{core}\left( \nu \right) $, it follows that $\mathcal{I}%
\not=\emptyset $ and $\mathrm{core}\left( \nu \right) \subseteq \mathcal{PI}$%
.\hfill $\square $

\smallskip

By the previous claim and Theorem \ref{thm:erg-PI}, the main statement
follows.

\smallskip

Finally, assume that $\nu $ is further ergodic. By Lemma \ref{lem:erg} and
since $f^{\star }\in B\left( \Omega ,\mathcal{G}\right) $ and $\nu $ is an
ergodic lower probability, it follows that%
\begin{equation*}
\nu \left( \left \{ \omega \in \Omega :\int_{\Omega }f^{\star }d\nu \leq
f^{\star }\left( \omega \right) \leq \int_{\Omega }f^{\star }d\bar{\nu}%
\right \} \right) =1.
\end{equation*}%
Since $\nu \left( \left \{ \omega \in \Omega :f^{\star }\left( \omega
\right) =\lim_{n}\frac{1}{n}\dsum \limits_{k=1}^{n}f\left( \tau ^{k-1}\left(
\omega \right) \right) \right \} \right) =1$ and $\nu $ is a lower
probability, this implies that%
\begin{equation*}
\nu \left( \left \{ \omega \in \Omega :\int_{\Omega }f^{\star }d\nu \leq
\lim_{n}\frac{1}{n}\dsum \limits_{k=1}^{n}f\left( \tau ^{k-1}\left( \omega
\right) \right) \leq \int_{\Omega }f^{\star }d\bar{\nu}\right \} \right) =1,
\end{equation*}%
proving the statement.\hfill $\blacksquare $

\smallskip

\noindent \textbf{Proof of Corollary \ref{cor:w-inv-erg}. }It is the proof
of the claim contained in the proof of Theorem \ref{thm:w-inv-erg}.\hfill $%
\blacksquare$

\smallskip

We next proceed by proving Theorem \ref{thm:sup-erg} and obtaining Corollary %
\ref{cor:erg} as a corollary of this former result. It is also possible to
provide a proof of Corollary \ref{cor:erg} as a consequence of Theorem \ref%
{thm:w-inv-erg}. By Theorem \ref{thm:w-inv-erg}, the extra assumption of $%
\left( \Omega ,\mathcal{F}\right) $ being standard yields the extra property
that $f^{\star }$ can be chosen to be the regular conditional expectation of 
$f$. Convexity and strong invariance imply that $\mathrm{core}\left( \nu
\right) \subseteq \mathcal{I}$. This yields that\ $\int_{\Omega }f^{\star
}d\nu =\int_{\Omega }fd\nu $ as well as $\int_{\Omega }f^{\star }d\bar{\nu}%
=\int_{\Omega }fd\bar{\nu}$. This, in turn, yields\ a sharper result under
the assumption of $\nu $ being ergodic.

\begin{lemma}
\label{lem:sup-sub}Let $\left \{ S_{n}\right \} _{n\in 
\mathbb{N}
}$ be a superadditive (resp., subadditive)\ sequence that satisfies (\ref%
{eq:bdd}) and $\mathcal{M}$ a compact subset of invariant probability
measures. If $\left \{ a_{n}\right \} _{n\in 
\mathbb{N}
}$ in $%
\mathbb{R}
$ is defined by $a_{n}=-\min_{P\in \mathcal{M}}\int_{\Omega }S_{n}dP$
(resp., $a_{n}=\max_{P\in \mathcal{M}}\int_{\Omega }S_{n}dP$) for all $n\in 
\mathbb{N}
$, then $\left \{ a_{n}\right \} _{n\in 
\mathbb{N}
}$ is subadditive, that is, $a_{n+k}\leq a_{n}+a_{k}$ for all $n,k\in 
\mathbb{N}
$.
\end{lemma}

\noindent \textbf{Proof. }Since $\left \{ S_{n}\right \} _{n\in 
\mathbb{N}
}$ satisfies (\ref{eq:bdd}), $\left \{ S_{n}\right \} _{n\in 
\mathbb{N}
}\subseteq B\left( \Omega ,\mathcal{F}\right) $. We just prove the statement
for the superadditive case, being the subadditive one similarly proven. If $%
\left \{ S_{n}\right \} _{n\in 
\mathbb{N}
}$ is superadditive and $\mathcal{M}$ is a compact subset of invariant
probability measures, then we have that $-a_{n+k}=\min_{P\in \mathcal{M}%
}\int_{\Omega }S_{n+k}dP\geq \min_{P\in \mathcal{M}}\int_{\Omega
}S_{n}+S_{k}\circ \tau ^{n}dP\geq \min_{P\in \mathcal{M}}\int_{\Omega
}S_{n}dP+\min_{P\in \mathcal{M}}\int_{\Omega }S_{k}\circ \tau
^{n}dP=\min_{P\in \mathcal{M}}\int_{\Omega }S_{n}dP+\min_{P\in \mathcal{M}%
}\int_{\Omega }S_{k}dP=-a_{n}-a_{k}$ for all $n,k\in 
\mathbb{N}
$, proving the statement.\hfill $\blacksquare $

\smallskip

\noindent \textbf{Proof of Theorem \ref{thm:sup-erg}. }Since $\nu $ is a
functionally invariant lower probability, we have that $\mathcal{M}\subseteq 
\mathcal{I}$. Define $\left \{ f_{n}\right \} _{n\in 
\mathbb{N}
}\subseteq B\left( \Omega ,\mathcal{F}\right) $ by $f_{n}=S_{n}/n$ for all $%
n\in 
\mathbb{N}
$. It follows that $\hat{f}_{n}\in B\left( \Omega ,\mathcal{G}\right) $ for
all $n\in 
\mathbb{N}
$. Since $\left \{ S_{n}\right \} _{n\in 
\mathbb{N}
}$ satisfies (\ref{eq:bdd}), it follows that there exists $\lambda \in 
\mathbb{R}
$ such that $-\lambda \leq f_{n},\hat{f}_{n}\leq \lambda $ for all $n\in 
\mathbb{N}
$. Define $f^{\star }\in B\left( \Omega ,\mathcal{G}\right) $ by $f^{\star
}=\sup_{n\in 
\mathbb{N}
}\hat{f}_{n}$ (resp., $f^{\star }=\inf_{n\in 
\mathbb{N}
}\hat{f}_{n}$). By Kingman's Subadditive Ergodic Theorem (see Dudley \cite[%
Theorem 10.7.1]{Dud} and \cite[Theorem 8.4]{Gra}) and since $W=\Omega $, we
have that $f^{\star }=\lim_{n}\hat{f}_{n}$ and $P\left( \left \{ \omega \in
\Omega :\lim_{n}\frac{S_{n}\left( \omega \right) }{n}=f^{\star }\left(
\omega \right) \right \} \right) =1$ for all $P\in \mathcal{M}$. Since $\nu $
is a lower probability, it follows that $\nu \left( \left \{ \omega \in
\Omega :\lim_{n}\frac{S_{n}\left( \omega \right) }{n}=f^{\star }\left(
\omega \right) \right \} \right) =1$, proving the main part of the statement.

\smallskip

1. If $\nu $ is convex and strongly invariant, then we have that $\mathrm{%
core}\left( \nu \right) \subseteq \mathcal{I}$ and%
\begin{equation}
\int_{\Omega }fd\nu =\min_{P\in \mathrm{core}\left( \nu \right)
}\int_{\Omega }fdP\qquad \forall f\in B\left( \Omega ,\mathcal{F}\right) .
\label{eq:cvx}
\end{equation}%
Consider the sequence $\left \{ a_{n}\right \} _{n\in 
\mathbb{N}
}$ defined by $a_{n}=-\int_{\Omega }S_{n}d\nu $ for all $n\in 
\mathbb{N}
$. By (\ref{eq:cvx}) and Lemma \ref{lem:sup-sub}, we have that $\left \{
a_{n}\right \} _{n\in 
\mathbb{N}
}$ is subadditive. It follows that (see \cite[Lemma 8.3]{Gra}) $\lim_{n}%
\frac{a_{n}}{n}=\inf_{n\in 
\mathbb{N}
}\frac{a_{n}}{n}$, that is,%
\begin{equation}
\lim_{n}\frac{-a_{n}}{n}=\sup_{n\in 
\mathbb{N}
}\frac{-a_{n}}{n}.  \label{eq:HW}
\end{equation}%
Recall that $\left \{ \hat{f}_{n}\right \} _{n\in 
\mathbb{N}
}$ is uniformly bounded. By Cerreia-Vioglio, Maccheroni, Marinacci, and
Montrucchio \cite[Theorem 22]{CMMM4}, (\ref{eq:HW}), and the main part of
the statement and since $\mathrm{core}\left( \nu \right) \subseteq \mathcal{I%
}$, we have that%
\begin{align*}
\int_{\Omega }f^{\star }d\nu & =\int_{\Omega }\lim_{n}\hat{f}_{n}d\nu
=\lim_{n}\int_{\Omega }\hat{f}_{n}d\nu =\lim_{n}\left[ \min_{P\in \mathrm{%
core}\left( \nu \right) }\int_{\Omega }\hat{f}_{n}dP\right] \\
& =\lim_{n}\left[ \min_{P\in \mathrm{core}\left( \nu \right) }\int_{\Omega
}f_{n}dP\right] =\lim_{n}\int_{\Omega }f_{n}d\nu =\lim_{n}\frac{\int_{\Omega
}S_{n}d\nu }{n} \\
& =\lim_{n}\frac{-a_{n}}{n}=\sup_{n\in 
\mathbb{N}
}\frac{-a_{n}}{n}=\sup_{n}\frac{\int_{\Omega }S_{n}d\nu }{n}=\sup_{n\in 
\mathbb{N}
}\int_{\Omega }f_{n}d\nu ,
\end{align*}%
proving point 1.

\smallskip

2. If $\nu $ is convex and strongly invariant, then we have that $\mathrm{%
core}\left( \nu \right) \subseteq \mathcal{I}$ and%
\begin{equation}
\int_{\Omega }fd\bar{\nu}=\max_{P\in \mathrm{core}\left( \nu \right)
}\int_{\Omega }fdP\qquad \forall f\in B\left( \Omega ,\mathcal{F}\right) .
\label{eq:ccv}
\end{equation}%
Consider the sequence $\left \{ a_{n}\right \} _{n\in 
\mathbb{N}
}$ defined by $a_{n}=\int_{\Omega }S_{n}d\bar{\nu}$. By (\ref{eq:ccv}) and
Lemma \ref{lem:sup-sub}, we have that $\left \{ a_{n}\right \} _{n\in 
\mathbb{N}
}$ is subadditive. It follows that (see \cite[Lemma 8.3]{Gra})%
\begin{equation}
\lim_{n}\frac{a_{n}}{n}=\inf_{n}\frac{a_{n}}{n}.  \label{eq:HW-ori}
\end{equation}%
Recall that $\left \{ \hat{f}_{n}\right \} _{n\in 
\mathbb{N}
}$ is uniformly bounded. By \cite[Theorem 22]{CMMM4}, (\ref{eq:HW-ori}), and
the main part of the statement and since $\mathrm{core}\left( \nu \right)
\subseteq \mathcal{I}$, we have that%
\begin{align*}
\int_{\Omega }f^{\star }d\bar{\nu}& =\int_{\Omega }\lim_{n}\hat{f}_{n}d\bar{%
\nu}=\lim_{n}\int_{\Omega }\hat{f}_{n}d\bar{\nu}=\lim_{n}\left[ \max_{P\in 
\mathrm{core}\left( \nu \right) }\int_{\Omega }\hat{f}_{n}dP\right] \\
& =\lim_{n}\left[ \max_{P\in \mathrm{core}\left( \nu \right) }\int_{\Omega
}f_{n}dP\right] =\lim_{n}\int_{\Omega }f_{n}d\bar{\nu}=\lim_{n}\frac{%
\int_{\Omega }S_{n}d\bar{\nu}}{n} \\
& =\lim_{n}\frac{a_{n}}{n}=\inf_{n}\frac{a_{n}}{n}=\inf_{n}\frac{%
\int_{\Omega }S_{n}d\bar{\nu}}{n}=\inf_{n\in 
\mathbb{N}
}\int_{\Omega }f_{n}d\bar{\nu},
\end{align*}%
proving point 2.

\smallskip

3. By Lemma \ref{lem:erg} and since $\nu $ is ergodic, it follows that%
\begin{equation*}
\nu \left( \left \{ \omega \in \Omega :\int_{\Omega }f^{\star }d\nu \leq
f^{\star }\left( \omega \right) \leq \int_{\Omega }f^{\star }d\bar{\nu}%
\right \} \right) =1.
\end{equation*}%
By the initial part of the proof, we have that $\nu \left( \left \{ \omega
\in \Omega :f^{\star }\left( \omega \right) =\lim_{n}\frac{S_{n}\left(
\omega \right) }{n}\right \} \right) =1$. Since $\nu $ is a lower
probability, this implies that%
\begin{equation*}
\nu \left( \left \{ \omega \in \Omega :\int_{\Omega }f^{\star }d\nu \leq
\lim_{n}\frac{S_{n}\left( \omega \right) }{n}\leq \int_{\Omega }f^{\star }d%
\bar{\nu}\right \} \right) =1,
\end{equation*}%
proving the statement.\hfill $\blacksquare $

\smallskip

\noindent \textbf{Proof of Corollary \ref{cor:erg}. }Pick $f\in B\left(
\Omega ,\mathcal{F}\right) $. It is immediate to see that $\left \{
S_{n}\right \} _{n\in 
\mathbb{N}
}$, defined by $S_{n}=\sum_{k=1}^{n}f\circ \tau ^{k-1}$ for all $n\in 
\mathbb{N}
$, is an additive sequence which satisfies (\ref{eq:bdd}). Since $\nu $ is
convex, continuous at $\Omega $,\ and strongly invariant, it is a
functionally invariant lower probability.\ Define $\left \{ f_{n}\right \}
_{n\in 
\mathbb{N}
}$ by $f_{n}=S_{n}/n$ for all $n\in 
\mathbb{N}
$. Note that $\hat{f}_{n}=\hat{f}$ for all $n\in 
\mathbb{N}
$. By the proof of Theorem \ref{thm:sup-erg},\textbf{\ }we have that $%
\lim_{n}\frac{S_{n}}{n}=\lim_{n}\hat{f}_{n}=\hat{f}$,\ $\nu -a.s.$, proving
the main statement and point 1 where $f^{\star }=\hat{f}$.

\smallskip

2. Since $\nu $ is convex and strongly invariant, then we have that $\mathrm{%
core}\left( \nu \right) \subseteq \mathcal{I}$ and $\int_{\Omega }fd\nu
=\min_{P\in \mathrm{core}\left( \nu \right) }\int_{\Omega }fdP$. By point 1
and since $\mathrm{core}\left( \nu \right) \subseteq \mathcal{I}$, we have
that $\int_{\Omega }fd\nu =\min_{P\in \mathrm{core}\left( \nu \right)
}\int_{\Omega }fdP=\min_{P\in \mathrm{core}\left( \nu \right) }\int_{\Omega }%
\hat{f}dP=\int_{\Omega }\hat{f}d\nu $, proving point 2. Note also that $%
\int_{\Omega }fd\bar{\nu}=\max_{P\in \mathrm{core}\left( \nu \right)
}\int_{\Omega }fdP=\max_{P\in \mathrm{core}\left( \nu \right) }\int_{\Omega }%
\hat{f}dP=\int_{\Omega }\hat{f}d\bar{\nu}$.

\smallskip

3. By point 3 of Theorem \ref{thm:sup-erg}\ and the proof of point 2, the
statement follows.\hfill $\blacksquare $

\smallskip

\noindent \textbf{Proof of Lemma \ref{lem:law}. }Consider a convex capacity $%
\nu $ and a process $\mathbf{f}$. It is immediate to see that $\nu _{\mathbf{%
f}}$ is a convex capacity. Next, consider $\left\{ C_{n}\right\} _{n\in 
\mathbb{N}
}\subseteq \sigma \left( \mathcal{C}\right) $ such that $C_{n}\uparrow 
\mathbb{R}
^{%
\mathbb{N}
}$. It follows that the sequence $\left\{ A_{n}\right\} _{n\in 
\mathbb{N}
}$, defined by $A_{n}=\mathbf{f}^{-1}\left( C_{n}\right) $ for all $n\in 
\mathbb{N}
$, is such that $A_{n}\uparrow \Omega $. Since $\nu $ is continuous at $%
\Omega $, we have that $\lim_{n}\nu _{\mathbf{f}}\left( C_{n}\right)
=\lim_{n}\nu \left( \mathbf{f}^{-1}\left( C_{n}\right) \right) =\lim_{n}\nu
\left( A_{n}\right) =1$, proving that $\nu _{\mathbf{f}}$ is continuous at $%
\mathbb{R}
^{%
\mathbb{N}
}$. Next, consider $C\in \mathcal{C}$. Then, there exist $k\in 
\mathbb{N}
$ and $E\in \mathcal{B}\left( 
\mathbb{R}
^{k}\right) $ such that $C=\left\{ x\in 
\mathbb{R}
^{%
\mathbb{N}
}:\left( x_{1},...,x_{k}\right) \in E\right\} $. Note that $\tau ^{-1}\left(
C\right) =\{x\in 
\mathbb{R}
^{%
\mathbb{N}
}:\left( x_{1},x_{2},...,x_{k+1}\right) \in 
\mathbb{R}
\times E\}$. Since $\mathbf{f}$ is stationary, it follows that%
\begin{align*}
\nu _{\mathbf{f}}\left( C\right) & =\nu \left( \mathbf{f}^{-1}\left(
C\right) \right) =\nu \left( \left\{ \omega \in \Omega :\left( f_{1}\left(
\omega \right) ,...,f_{k}\left( \omega \right) \right) \in E\right\} \right) 
\\
& =\nu \left( \left\{ \omega \in \Omega :\left( f_{2}\left( \omega \right)
,...,f_{k+1}\left( \omega \right) \right) \in E\right\} \right)  \\
& =\nu \left( \left\{ \omega \in \Omega :\left( f_{1}\left( \omega \right)
,f_{2}\left( \omega \right) ,...,f_{k+1}\left( \omega \right) \right) \in 
\mathbb{R}
\times E\right\} \right)  \\
& =\nu \left( \mathbf{f}^{-1}\left( \tau ^{-1}\left( C\right) \right)
\right) =\nu _{\mathbf{f}}\left( \tau ^{-1}\left( C\right) \right) .
\end{align*}%
Since $C\in \mathcal{C}$ was arbitrarily chosen, it follows that $\mathcal{C}%
\subseteq \{C\in \sigma \left( \mathcal{C}\right) :\nu _{\mathbf{f}}\left(
C\right) =\nu _{\mathbf{f}}\left( \tau ^{-1}\left( C\right) \right)
\}\subseteq \sigma \left( \mathcal{C}\right) $. Since $\nu _{\mathbf{f}}$ is
convex and continuous at $%
\mathbb{R}
^{%
\mathbb{N}
}$, we have that $\{C\in \sigma \left( \mathcal{C}\right) :\nu _{\mathbf{f}%
}\left( C\right) =\nu _{\mathbf{f}}\left( \tau ^{-1}\left( C\right) \right)
\}$ is a monotone class. By the Monotone Class Theorem (see \cite[Theorem 3.4%
]{Bil}), it follows that $\sigma \left( \mathcal{C}\right) =\left\{ C\in
\sigma \left( \mathcal{C}\right) :\nu _{\mathbf{f}}\left( C\right) =\nu _{%
\mathbf{f}}\left( \tau ^{-1}\left( C\right) \right) \right\} $, that is, $%
\nu _{\mathbf{f}}$ is shift invariant. Define $\mathcal{H}%
=\dbigcap\limits_{k=1}^{\infty }\sigma \left( \mathcal{C}_{k+1}^{\infty
}\right) \cap \sigma \left( \mathcal{C}\right) $.\footnote{$\mathcal{C}%
_{k+1}^{\infty }$ is the class of cylinders such that%
\begin{equation*}
C=\left\{ x\in 
\mathbb{R}
^{%
\mathbb{N}
}:\left( x_{1},...,x_{k},x_{k+1},...,x_{k^{\prime }}\right) \in 
\mathbb{R}
^{k}\times E\right\} 
\end{equation*}%
where $k^{\prime }>k$ and $E\in \mathcal{B}(%
\mathbb{R}
^{k^{\prime }-k})$.} Note that $\mathbf{f}^{-1}\left( \mathcal{H}\right)
\subseteq \mathcal{T}$. Thus, $\nu _{\mathbf{f}}\left( \mathcal{H}\right)
=\left\{ 0,1\right\} $ if $\nu \left( \mathcal{T}\right) =\left\{
0,1\right\} $. Let $\mathcal{G}$ be the $\sigma $-algebra of shift invariant
events. It is well known that $\mathcal{G}\subseteq \mathcal{H}$. In light
of these observations, it is immediate to see that if $\nu \left( \mathcal{T}%
\right) =\left\{ 0,1\right\} $, then $\nu _{\mathbf{f}}\left( \mathcal{G}%
\right) =\left\{ 0,1\right\} $, that is, $\nu _{\mathbf{f}}\ $is
ergodic.\hfill $\blacksquare $

\smallskip

\noindent \textbf{Proof of Theorem \ref{thm:SLLN}. }By induction and since $%
\mathbf{f}$ is stationary, it follows that for each $k\in 
\mathbb{N}
$ and for each Borel subset $B$ of $%
\mathbb{R}
$%
\begin{equation}
\nu \left( \left \{ \omega \in \Omega :f_{1}\left( \omega \right) \in
B\right \} \right) =\nu \left( \left \{ \omega \in \Omega :f_{k}\left(
\omega \right) \in B\right \} \right) .  \label{eq:dib}
\end{equation}%
By (\ref{eq:dib}), this implies that for each $k\in 
\mathbb{N}
$ and for each Borel subset $B$ of $%
\mathbb{R}
$%
\begin{equation*}
\nu _{\mathbf{f}}\left( \left \{ x\in 
\mathbb{R}
^{%
\mathbb{N}
}:x_{k}\in B\right \} \right) =\nu \left( \left \{ \omega \in \Omega
:f_{k}\left( \omega \right) \in B\right \} \right) =\nu \left( \left \{
\omega \in \Omega :f_{1}\left( \omega \right) \in B\right \} \right) .
\end{equation*}%
In particular, since $\left \{ f_{n}\right \} _{n\in 
\mathbb{N}
}\subseteq B\left( \Omega ,\mathcal{F}\right) $, it follows that there
exists $m\in 
\mathbb{R}
$ such that $-m1_{\Omega }\leq f_{1}\leq m1_{\Omega }$. If we replace $B$
with $\left[ -m,m\right] $, then we can conclude that%
\begin{equation}
\nu _{\mathbf{f}}\left( \left \{ x\in 
\mathbb{R}
^{%
\mathbb{N}
}:x_{k}\in \left[ -m,m\right] \right \} \right) =\nu \left( \left \{ \omega
\in \Omega :f_{1}\left( \omega \right) \in \left[ -m,m\right] \right \}
\right) =1\quad \forall k\in 
\mathbb{N}
.  \label{eq:bdd-v}
\end{equation}%
Define $\pi :%
\mathbb{R}
^{%
\mathbb{N}
}\rightarrow 
\mathbb{R}
$ by%
\begin{equation*}
\pi \left( x\right) =\left \{ 
\begin{array}{cc}
x_{1} & if\text{ }x_{1}\in \left[ -m,m\right] \\ 
0 & otherwise%
\end{array}%
\right. \qquad \forall x\in 
\mathbb{R}
^{%
\mathbb{N}
}.
\end{equation*}%
It is immediate to see that $\pi \in B\left( 
\mathbb{R}
^{%
\mathbb{N}
},\sigma \left( \mathcal{C}\right) \right) $. Note also that%
\begin{equation}
\dbigcap \limits_{k=1}^{\infty }\left \{ x\in 
\mathbb{R}
^{%
\mathbb{N}
}:x_{k}\in \left[ -m,m\right] \right \} \subseteq \dbigcap
\limits_{n=1}^{\infty }\left \{ x\in 
\mathbb{R}
^{%
\mathbb{N}
}:\frac{1}{n}\sum_{k=1}^{n}\pi \left( \tau ^{k-1}\left( x\right) \right) =%
\frac{1}{n}\sum_{k=1}^{n}x_{k}\right \} .  \label{eq:inc}
\end{equation}%
By (\ref{eq:bdd-v}) and (\ref{eq:inc})\ and since $\nu _{\mathbf{f}}$ is a
convex capacity which is further continuous at $%
\mathbb{R}
^{%
\mathbb{N}
}$, it follows that%
\begin{equation}
\nu _{\mathbf{f}}\left( \dbigcap \limits_{n=1}^{\infty }\left \{ x\in 
\mathbb{R}
^{%
\mathbb{N}
}:\frac{1}{n}\sum_{k=1}^{n}\pi \left( \tau ^{k-1}\left( x\right) \right) =%
\frac{1}{n}\sum_{k=1}^{n}x_{k}\right \} \right) =1.  \label{1}
\end{equation}%
By Theorem \ref{thm:w-inv-erg} and since $\nu _{\mathbf{f}}$ is shift
invariant and ergodic, we have that there exists $\pi ^{\star }\in B\left( 
\mathbb{R}
^{%
\mathbb{N}
},\mathcal{G}\right) $ such that%
\begin{equation}
\nu _{\mathbf{f}}\left( \left \{ x\in 
\mathbb{R}
^{%
\mathbb{N}
}:\int_{%
\mathbb{R}
^{%
\mathbb{N}
}}\pi ^{\star }d\nu _{\mathbf{f}}\leq \lim_{n}\frac{1}{n}\sum_{k=1}^{n}\pi
\left( \tau ^{k-1}\left( x\right) \right) =\pi ^{\star }\left( x\right) \leq
\int_{%
\mathbb{R}
^{%
\mathbb{N}
}}\pi ^{\star }d\bar{\nu}_{\mathbf{f}}\right \} \right) =1.  \label{2}
\end{equation}%
By (\ref{1}) and (\ref{2}) and since $\nu _{\mathbf{f}}$ is convex, we can
conclude that%
\begin{equation}
\nu _{\mathbf{f}}\left( \left \{ x\in 
\mathbb{R}
^{%
\mathbb{N}
}:\int_{%
\mathbb{R}
^{%
\mathbb{N}
}}\pi ^{\star }d\nu _{\mathbf{f}}\leq \lim_{n}\frac{1}{n}%
\sum_{k=1}^{n}x_{k}=\pi ^{\star }\left( x\right) \leq \int_{%
\mathbb{R}
^{%
\mathbb{N}
}}\pi ^{\star }d\bar{\nu}_{\mathbf{f}}\right \} \right) =1.
\label{eq:dis-SLLN}
\end{equation}%
Let\ $E=\left \{ x\in 
\mathbb{R}
^{%
\mathbb{N}
}:\lim_{n}\frac{1}{n}\sum_{k=1}^{n}\pi \left( \tau ^{k-1}\left( x\right)
\right) =\pi ^{\star }\left( x\right) \right \} $ and $\pi _{n}=\frac{1}{n}%
\sum_{k=1}^{n}\pi \left( \tau ^{k-1}\right) $ for all $n\in 
\mathbb{N}
$. By (\ref{2}), we have that $P\left( E\right) =1$ for all $P\in \mathrm{%
core}\left( \nu _{\mathbf{f}}\right) $. By construction, $\left \{ 1_{E}\pi
_{n}\right \} _{n\in 
\mathbb{N}
}\subseteq B\left( 
\mathbb{R}
^{%
\mathbb{N}
},\sigma \left( \mathcal{C}\right) \right) $ is a uniformly bounded sequence
which converges pointwise to $1_{E}\pi ^{\star }$. By \cite[Theorem 22]%
{CMMM4} and since $\nu _{\mathbf{f}}$ is convex and $P\left( E\right) =1$
for all $P\in \mathrm{core}\left( \nu _{\mathbf{f}}\right) $, this implies
that%
\begin{equation}
\int_{%
\mathbb{R}
^{%
\mathbb{N}
}}\pi ^{\star }d\nu _{\mathbf{f}}=\int_{%
\mathbb{R}
^{%
\mathbb{N}
}}1_{E}\pi ^{\star }d\nu _{\mathbf{f}}=\int_{%
\mathbb{R}
^{%
\mathbb{N}
}}\lim_{n}1_{E}\pi _{n}d\nu _{\mathbf{f}}=\lim_{n}\int_{%
\mathbb{R}
^{%
\mathbb{N}
}}1_{E}\pi _{n}d\nu _{\mathbf{f}}=\lim_{n}\int_{%
\mathbb{R}
^{%
\mathbb{N}
}}\pi _{n}d\nu _{\mathbf{f}}.  \label{eq:lim}
\end{equation}%
Next, since $\nu _{\mathbf{f}}$ is convex and shift invariant, note that for
each $n\in 
\mathbb{N}
$%
\begin{equation*}
\int_{%
\mathbb{R}
^{%
\mathbb{N}
}}\pi _{n}d\nu _{\mathbf{f}}=\int_{%
\mathbb{R}
^{%
\mathbb{N}
}}\frac{1}{n}\sum_{k=1}^{n}\pi \left( \tau ^{k-1}\right) d\nu _{\mathbf{f}%
}\geq \frac{1}{n}\sum_{k=1}^{n}\int_{%
\mathbb{R}
^{%
\mathbb{N}
}}\pi \left( \tau ^{k-1}\right) d\nu _{\mathbf{f}}=\int_{%
\mathbb{R}
^{%
\mathbb{N}
}}\pi d\nu _{\mathbf{f}}.
\end{equation*}%
By (\ref{eq:lim}), it follows that $\int_{%
\mathbb{R}
^{%
\mathbb{N}
}}\pi ^{\star }d\nu _{\mathbf{f}}\geq \int_{%
\mathbb{R}
^{%
\mathbb{N}
}}\pi d\nu _{\mathbf{f}}$. A similar argument yields that $\int_{%
\mathbb{R}
^{%
\mathbb{N}
}}\pi ^{\star }d\bar{\nu}_{\mathbf{f}}\leq \int_{%
\mathbb{R}
^{%
\mathbb{N}
}}\pi d\bar{\nu}_{\mathbf{f}}$. Finally, since $\int_{%
\mathbb{R}
^{%
\mathbb{N}
}}\pi d\nu _{\mathbf{f}}=\int_{\Omega }f_{1}d\nu $ and $\int_{%
\mathbb{R}
^{%
\mathbb{N}
}}\pi d\bar{\nu}_{\mathbf{f}}=\int_{\Omega }f_{1}d\bar{\nu}$, by (\ref%
{eq:dis-SLLN}), we can conclude that%
\begin{align*}
1& =\nu _{\mathbf{f}}\left( \left \{ x\in 
\mathbb{R}
^{%
\mathbb{N}
}:\int_{%
\mathbb{R}
^{%
\mathbb{N}
}}\pi d\nu _{\mathbf{f}}\leq \lim_{n}\frac{1}{n}\sum_{k=1}^{n}x_{k}\leq
\int_{%
\mathbb{R}
^{%
\mathbb{N}
}}\pi d\bar{\nu}_{\mathbf{f}}\right \} \right) \\
& =\nu \left( \left \{ \omega \in \Omega :\int_{\Omega }f_{1}d\nu \leq
\lim_{n}\frac{1}{n}\sum_{k=1}^{n}f_{k}\left( \omega \right) \leq
\int_{\Omega }f_{1}d\bar{\nu}\right \} \right) ,
\end{align*}%
proving the statement.\hfill $\blacksquare $

\end{document}